\documentclass[12pt]{amsart}

\usepackage{amssymb,amsmath,amsthm,amsfonts,times,hyperref,mathrsfs,multirow}

\usepackage{epsfig}
\usepackage{amsmath}
\usepackage{amssymb}

\begin{document}
\def\eqn#1{(\ref{eq:#1})}
\def\ssn{\vskip 2pt\noindent}
\def\sn{\vskip 10pt\noindent}
\def\bn{\vskip 20pt\noindent}
\def\w{\omega}
\def\qed{\penalty-250%
\hbox to 6pt{\hfill}\hfill\llap{%
\vbox{\hrule\hbox{\vrule height6pt\hskip6pt\vrule}\hrule}}}
\def\tentt{\tt}
\newcommand{\qbin}[2]{\left[\begin{matrix} #1 \\ #2 \end{matrix} \right]}
\newcommand{\bin}[2]{\left(\begin{matrix} #1 \\ #2 \end{matrix} \right)}

\newcommand\mybibitem[1]{\bibitem{#1}}
\renewcommand{\MR}[1]{\href{http://www.ams.org/mathscinet-getitem?mr={#1}}{MR{#1}}}
\numberwithin{equation}{section}

\newtheorem{thm}{Theorem}

%
%

\title[Basis partitions]
{Basis partitions and their signature}

\author{Krishnaswami Alladi}
\address{Department of Mathematics, University of Florida, Gainesville,
FL 32611-8105}
\email{alladik@ufl.edu}
                                   

\dedicatory{Dedicated to Steve Milne on the occasion of his 75th birthday}

\date{\today}                   
\thanks{Research supported in part by NSA Grants MSPF-06G-150 and 
MSPF-08G-154}

\begin{abstract}
Basis partitions are minimal partitions 
corresponding to successive rank vectors. We show combinatorially 
how basis partitions can be generated from primary partitions 
which are equivalent to the Rogers-Ramanujan partitions. This 
leads to the definition of a signature of a basis partition that we 
use to explain certain parity results. We then study a special class of 
basis partitions which we term as complete. Finally we discuss basis 
partitions and minimal basis partitions among partitions with non-repeating 
odd parts by representing them using 2-modular graphs.
\end{abstract}

\maketitle

\section{Introduction}

For a partition $\pi: b_1\ge b_2\ge ...\ge b_\nu$, its Ferrers graph 
is an array of left-justified rows with $b_i$ nodes in the i-th row, 
$i=1,2,...,\nu$. 
The conjugate $\pi^*$ of $\pi$ is 
obtained 
by counting the nodes in the columns of the Ferrers graph of $\pi$. 
Let 
$\pi^*$ be $c_1\ge c_2 ... \ge c_{b_1}$.  

In every Ferrers graph, there is a largest square of nodes called the 
{\it{Durfee square}} which is the same for the conjugate as well. Given 
a partition $\pi$ as above, let it have a $k\times k$ Durfee square. 
The successive rank vector of $\pi$ is defined as {\bf{r}}= $(r_1, r_2, ..., 
r_k)$, where 
\begin{equation}
r_i=b_i-c_i, \quad \text{for} \quad i=1,2,...k. \label{eq:1.1}
\end{equation}
 
Note that $r_1=b_1-c_1$ is the {\it{rank}} that Dyson \cite{Dy44} used 
to combinatorially explain two of Ramanujan's congruences. This led
to the terminology {\it{successive ranks}} for the differences $b_i-c_i$ 
which are equal to hook differences. Their study has led to the discovery of 
important partition identities, the results in Andrews et al \cite{An-Ba-etal0087} being a good 
example. 

Given {\bf{r}}, there are, in general, several partitions $\pi$ 
which will have {\bf{r}} as the successive rank vector. Among all 
partitions $\pi$ having a prescribed successive rank vector, consider 
the partition for which the sum of the parts $\sigma(\pi)=b_1+b_2+...+
b_{\nu}$ (the integer being partitioned) is minimal. Gupta \cite{Gu0076} called 
such minimal partitions as {\it{basis partitions}} and proved their 
existence. He also noticed that if $\pi$ is the basis partition associated 
with a rank vector {\bf{r}}= $(r_1, r_2, ..., r_k)$ and if $n>\sigma(\pi)$, then the number 
of partitions of $n$ which have the same rank vector {\bf{r}} is equal to 
the number of partitions of $(n-\sigma(\pi))/2$ into $\le k$ parts. 

We provide here a combinatorial procedure to generate basis 
partitions from primary partitions which are equivalent to the 
Rogers-Ramanujan partitions (see Section 3), and use that for a new study 
of basis 
partitions among partitions with non-repeating odd parts (Section 9). 
Hirschhorn \cite{Hi0099} had noticed 
that the number of basis partitions of an integer $n$ can be expressed 
as a weighted count over the Rogers-Ramanujan partitions of $n$, where 
the weights are powers of 2. He deduced this by suitable interpretation 
of the generating function of basis partitions obtained by Nolan, Savage, 
and Wilf \cite{No-Sa-Wi0098}. 


Our approach provides a very natural combinatorial proof of 
Hirschhorn's weighted identity 
as well as a refinement of the generating function of basis 
partitions that leads to a statistic which we call the {\it{signature}}.
A study of the parity of the signature yields a partial theta series involving 
the squares (Section 5). This parity result has nice 
counterparts for partitions with non-repeating odd parts (Section 9).

\section{The generating function of basis partitions}

\medskip
The Ferrers graph of any partition $\pi$ may be split into three 
parts as follows: (i) the Durfee square $D(\pi)$, (ii) the partition 
$\pi_r$ which is the part to the right of the Durfee square, and 
(iii) the partition $\pi_b$ which is the portion below the Durfee 
square. If in a partition there is a column of $\pi_r$ which is 
equal in length to a row of $\pi_b$, then we may remove this 
row and column to get a smaller partition which will have the 
same successive rank vector. Thus a basis partition, by definition 
of its minimality, will not have such a redundancy. Thus we have 
the following {\it{Characterization: $\pi$ is a basis partition 
if and only if no column of $\pi_r$ equals any row of $\pi_b$.}} 
Conversely, given basis partition with rank vector {\bf{r}} 
of length $k$, all partitions having the same rank vector can be 
generated by inserting columns of length $\le k$ to the right of 
the Durfee square, and inserting rows below the Durfee square 
exactly equal in length and multiplicity to the inserted columns. 
This observation can be used to compute the generating function 
of basis partitions. 

Let $b(n;k)$ denote the number of basis partitions whose rank vector 
is of length $k$. Let $p(n;k)$ denote the number of partitions 
having a $k\times k$ Durfee square. It is well known that
\begin{equation}
\sum_{n}p(n;k)q^n=\frac{q^{k^2}}{\{(1-q)(1-q^2)...(1-q^k)\}^2}.\label{eq:2.1}
\end{equation}
Also, from the observation in the preceding paragraph, we have 
\begin{equation}
\sum_{n}p(n;k)q^n=\frac{1}{(1-q^2)(1-q^4)...(1-q^{2k})}
\sum_{n}b(n;k)q^n.\label{eq:2.2}
\end{equation}
From \eqn{2.1} and \eqn{2.2} it follows that 
\begin{equation}
\sum_{n}b(n;k)q^n=\frac{q^{k^2}(1+q)(1+q^2)...(1+q^k)}
{(1-q)(1-q^2)...(1-q^k)}.\label{eq:2.3}
\end{equation}


The derivation of of \eqn{2.3} given above is due to Hirschhorn \cite{Hi0099}. 
In \cite{No-Sa-Wi0098}, the above characterization of basis partitions is observed, 
but the derivation of \eqn{2.3} is by other means. Andrews \cite{An2015} has studied basis partitions and a certain polynomial assosiated with \eqn{2.3}, but he interpreted
\eqn{2.3} in terms of overpartitions. 
 
Rogers-Ramanujan partitions $\pi: b_1+b_2+...+b_k$ are those whose 
parts differ by at least 2. That is 
\begin{equation}
b_i-b_{i+1}\ge 2, \quad \text{for} \quad i=1,2,..., k-1.\label{eq: 2.4}
\end{equation}
It is well known \cite{Al2010} that the generating function of the Rogers-Ramanujan 
partitions into exactly $k$ parts is
\begin{equation}
\frac{q^{k^2}}{(1-q)(1-q^2)...(1-q^k)}.\label{eq:2.5}
\end{equation}
The way one realizes this is to note that the minimal partition into 
$k$ parts differing by $\ge 2$ is $1+3+5+...+ (2k-1)=k^2$. This accounts 
for the term $q^{k^2}$ in the numerator in \eqn{2.5}. One can construct all 
partitions into $k$ parts that differ by $\ge 2$ from the minimal partition 
by inserting columns of length $\le k$ into the Ferrers graph. This accounts 
for the factors $(1-q)(1-q^2)...(1-q^k)$ in the denominator in \eqn{2.5}.      
     
Hirschhorn \cite{Gu0076} notes that since
\begin{equation}
\frac{(1+q^j)}{(1-q^j)}=1+2(q^j+q^{2j}+...+q^{nj}+...),\label{eq:2.6}
\end{equation}
the insertion of the columns of of a given length $j$ into the minimal 
partition increases the size of the gap and the increase is to be counted 
with weight 2. Thus he deduces

\medskip
\noindent
\begin{thm}[{ Hirschhorn \cite{Hi0099} }]
\label{thm:1}
Let {\bf{R}} {{denote the set of Rogers-Ramanujan partitions. For 
each partition}} $\pi\in${\bf{R}}, 
\begin{equation}
\pi: b_1+b_2+...+b_k, \quad b_i-b_{i+1}\ge 2, \quad \text{for} \quad 
i=1,2,..., k-1, \quad \text{and} \quad b_k\ge 1, \label{eq:2.7}
\end{equation}
{{let its weight be}} $\omega(\pi)=2^t$, {{where t is the number of 
strict inequalities in \eqn{2.7}. Then}}
$$
b(n,k)=\sum_{\pi\in R, \sigma(\pi)=n, \nu(\pi)=k}\omega(\pi).
$$
Here $\nu(\pi)$ denotes the number of parts of a partition $\pi$.   
\end{thm}

\bigskip
\section{Constructing basis partitions from primary partitions}

\medskip
We define a {\it{primary}} partition $\pi$ to be one for which $\pi_b$ is 
empty, that is, there are no nodes below the Durfee square. It is well known 
\cite{An-book} that the expression in \eqn{2.5} is the generating function for partitions 
into $k$ parts, each part $\ge k$. This is the same as saying that the 
expression in \eqn{2.5} is the generating function of primary partitions into 
$k$ parts. The primary partitions are bijectively equivalent to the 
Rogers-Ramanujan partitions as the following correspondence shows: Given 
the Ferrers graph of a primary partition $\pi$ into $k$ parts, consider 
the partition $\pi'$ obtained by counting nodes along the hooks of the 
graph. Then $\pi'$ is a Rogers-Ramanujan partition, that is, a partition 
into parts that differ by $\ge 2$. This correspondence can be reversed. 

The term {\it{primary partitions}} is ours \cite{Al0097}. Although it is well known that 
the expression in \eqn{2.5} is the generating function of these partitions, 
the most common common interpretation and use of \eqn{2.5} has been 
its role as the generating function of the Rogers-Ramanujan partitions. 
We noticed \cite{Al0097}, \cite{Al0098}, that primary partitions play a crucial role in the 
theory of weighted partition identities. Indeed we established a variety 
of weighted partition identities by performing {\it{sliding operations}}
(defined below) on primary partitions. This is also what motivated us to 
coin the term primary partitions.

Note that since a primary partition $\pi$ has no nodes below the Durfee 
square, we have trivially that no column of $\pi_r$ equals a row of 
$\pi_b$. Thus every primary partition is automatically a basis partition.     
We will now show how to generate all basis partitions from primary 
partitions. This will provide a combinatorial proof of Theorem 1. 
Before providing this construction, we reformulate Theorem 1 as:

\medskip

\noindent
\begin{thm}
\label{thm:2} 
{{Let}} {\bf{P}} {{denote the set of primary partitions. For 
each partition}} $\pi\in${\bf{P}}, 
\begin{equation}
\pi: b_1+b_2+...+b_k, \quad b_i-b_{i+1}\ge 0, \quad \text{for} \quad 
i=1,2,..., k-1, \quad \text{and} \quad b_k\ge k, \label{eq:3.1}
\end{equation}
{{let its weight be}} $\omega(\pi)=2^t$, {{where t is the number of 
strict inequalities in \eqn{3.1}. Then}}
$$
b(n,k)=\sum_{\pi\in P, \sigma(\pi)=n, \nu(\pi)=k}\omega(\pi).
$$
\end{thm}
\medskip
{\bf{Proof:}} Given a primary partition $\pi$, suppose a column of 
$\pi_r$ is moved and placed as a row below the Durfee square. We call 
this a {\it{sliding operation}}. On a given primary partition, several 
sliding operations can be performed. The following are invariants under 
the sliding operation: (i) $\sigma(\pi)$, the number being partitioned, 
(ii) $D(\pi)$, the Durfee square of the partition, and (iii) the hook 
lengths of the partition. This last invariant is important here 
because the underlying Rogers-Ramanujan partition obtained by counting 
nodes along hooks remains invariant under sliding operations. 

Suppose we consider columns of a certain length $\ell$ in the portion 
$\pi_r$ of the primary partition $\pi$. If a column of length $\ell$ 
is moved down and placed as a row below the Durfee square, then we will 
not have have a basis partition if there are columns of length $\ell$
remaining in the portion to the right of the Durfee square, because 
this violates the characterization of basis partitions (see Section 2). 
But then, if we move {\it{all}} columns of length $\ell$ and place them as 
rows below the Durfee square, then the characterization condition is 
confirmed and we have a basis partition. Thus under the sliding operation, 
we have a choice of either not to move any column of length $\ell$ or 
move {\it{all}} columns of length $\ell$ to form basis partitions. Thus we 
have two choices for every block of columns of a given length in $\pi_r$. 
Finally we note that under the sliding operation, movement of columns 
of a certain length is {\it{independent}} of the movement of columns 
of a different length. Thus if the primary partition $\pi$ contains 
$t$ sets of columns of different lengths in $\pi_r$, then it spawns
$2^t$ basis partitions (including itself) under the sliding operation. 
All basis partitions can be obtained from primary partitions this way. 
The number $t$ of sets of columns of different lengths in 
$\pi_r$ corresponds to the number of strict inequalities in \eqn{3.1} 
(and in \eqn{2.7}) and so this combinatorial construction provides a 
proof of Theorem 2 (and of Theorem 1). 

\bigskip

\section{Constructing a basis partition of a successive rank vector}

\medskip
We will now discuss how to construct a basis partition corresponding to 
a successive rank vector {\bf{r}}=$(r_1,r_2,...,r_k)$. The construction 
will be illustrated with the example 
$$
\mathbf{r}=(r_1,r_2,r_3,r_4,r_5)=(3,2,-1,4,-3).
$$

\noindent
{\bf{Step 1:}} First form a $k\times k$ Durfee square.  
\medskip

\begin{tabbing}
XXXXXXXXXXXXXXXXXXXXXXXX \= YYYYYYYYYYYYYYYYYY\kill
\vbox to 1truein{\hbox{Basis partition corresponding to}
                    \hbox{${\bold r}=(0,0,0,0,0)$}
                    \vfill} 
\>
 \hbox to 1truein{\epsfxsize=1truein\epsfbox{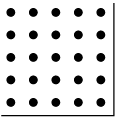}}\\
\end{tabbing}

\noindent
{\bf{Step 2:}} If $r_k>0$, form $r_k$ columns of length k to the right 
of the Durfee square. If $r_k<0$, form $|r_k|$ rows of length $k$ below 
the Durfee square. The resulting partition has all successive ranks equal 
to $r_k$. If $r_k=0$, then move on to $r_{k-1}$. 

\medskip
\begin{tabbing}
XXXXXXXXXXXXXXXXXXXXXXXX \= YYYYYYYYYYYYYYYYYY\kill
\vbox to 1.5truein{\hbox{$r_5=-3$} 
                    \hbox{Basis partition corresponding to}
                    \hbox{${\bold r}=(-3,-3,-3,-3,-3)$}
                    \vfill} 
\>
 \hbox to 1truein{\epsfxsize=1truein\epsfbox{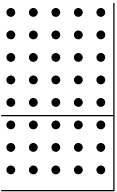}}\\
\end{tabbing}

\noindent
{\bf{Step 3:}} Next consider the difference $r_{k-1}-r_k=\delta_{k-1}$. 
If this difference is $>0$, place $\delta _{k-1}$ columns of length 
$k-1$ to the right of the Durfee square. If the difference is $<0$, 
place $|\delta_{k-1}|$ rows of length $k-1$ below the Durfee square. 
Note that all successive ranks up to the $k-1$-st are equal to 
$r_{k-1}$, and the last rank is $r_k$. If $\delta_{k-1}=0$, move on to 
$r_{k-2}$.          

\medskip

\begin{tabbing}
XXXXXXXXXXXXXXXXXXXXXXX \= YYYYYYYYYYYYYYYYYY\kill
\vbox to 1.5truein{\hbox{$r_4 - r_5=7 = \delta_4$}
                    \hbox{Basis partition corresponding to}
                    \hbox{${\bold r}=(4,4,4,4,-3)$}
                    \vfill} 
\>
 \hbox to 2.2truein{\epsfxsize=2.2truein\epsfbox{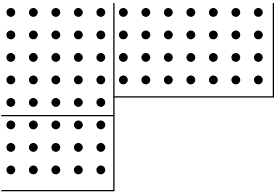}}\\
\end{tabbing}

\noindent
{\bf{Step 4 (General):}} Consider the difference $r_i-r_{i+1}=\delta_i$. 
If this difference is $>0$, place $\delta_i$ columns of length $i$ to the 
right of the Durfee square. If the difference is $<0$, place $|\delta_i|$ 
rows of length $i$ below the Durfee square. If $\delta_i=0$, proceed 
to $r_{i-1}$. Note that all successive ranks up to the i-th are equal 
to $r_i$, but from then on the successive ranks are $r_{i+1}, ..., r_k$. 

\noindent
{\bf{Step 5:}} Complete the construction by considering the final case 
$j=1$ in Step 4. We end up with a basis partition with successive rank 
vector $(r_1, r_2, ..., r_k)$ as desired.

\medskip

\begin{tabbing}
XXXXXXXXXXXXXXXXXXXXXXX \= YYYYYYYYYYYYYYYYYY\kill
\vbox to 2.9truein{\hbox{$\delta_1=5$, $\delta_2=-3$,}
                      \hbox{$\delta_3=-3$, $\delta_4=7$, $\delta_5=-3$}
                    \hbox{Basis partition corresponding to}
                    \hbox{${\bold r}=(3,-2,1,4,-3)$}
                    \vfill}
\>
 \hbox to 3.5truein{\epsfxsize=3.5truein\epsfbox{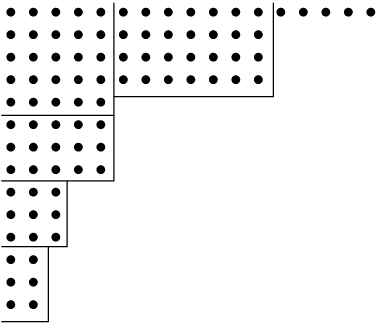}}\\
\end{tabbing}


\bigskip

\section{The signature of basis partitions and a theta series}

\medskip
The combinatorial proof given in Section 3, and the construction of 
the basis partition corresponding to successive rank vector given in 
Section 4, lead to the definition of the {\it{signature}} $\Sigma(\pi)$ 
of a basis partition to be the number of differences $r_i-r_{i+1}=\delta_i$ 
that are negative. Note that 
\begin{equation}
\Sigma (\pi)=\nu_d(\pi_b),\label{eq:5.1}
\end{equation}
where $\nu_d(\pi_b)$ is the number of different parts of $\pi_b$. The 
combinatorial constructions of Sections 3 and 4 enable us to keep track 
of this signature and obtain a refined generating function of basis 
partitions. More precisely, if $b(n,k,s)$ denotes the number of basis 
partitions of $n$ with a $k\times k$ Durfee square and with signature 
$s$, then for a fixed $k$, we get the following refinement of \eqn{2.3}:  
\begin{equation}
\sum_{n,s}b(n,k,s)q^nz^s=\frac{q^{k^2}(1+zq)(1+zq^2)...(1+zq^k)}
{(1-q)(1-q^2)...(1-q^k)}.\label{eq:5.2}
\end{equation}
Observe that $z=-1$ makes the right hand side of \eqn{5.2} collapse to 
$q^{k^2}$. Thus we get the partial theta series identity
\begin{equation}
\sum_{n,k,s}b(n,k,s)q^n(-1)^s=\sum^{\infty}_{k=0}q^{k^2}.\label{eq:5.3}
\end{equation}
This has the following partition interpretation:
\begin{thm}
\label{thm:3}
{{For an integer n, consider the difference between 
the number of basis partitions of n with even signature and the number of 
basis partitions with odd signature. This difference is 1 if n is 
a perfect square, and 0 otherwise.}}
\end{thm}
\medskip

Actually, the discussion of the signature leads to a refinement of 
Theorems 1 and 2, where the weight $2^t$ is replaced by $(1+z)^t$ so 
that Theorems 1 and 2 correspond to the case $z=1$ and Theorem 3 
follows as the special case $z=-1$. We state these refinements as 
Theorem 4 parts (i) and (ii) below.   

\begin{thm}
\label{thm:4} 
\leavevmode
\begin{enumerate}
\item[(i)]
{{Let}} {\bf{R}} {{denote the set of Rogers-Ramanujan partitions. For 
each partition}} $\pi\in${\bf{R}}, 
\begin{equation}
\pi: b_1+b_2+...+b_k, \quad b_i-b_{i+1}\ge 2, \quad \text{for} \quad 
i=1,2,..., k-1, \quad \text{and} \quad b_k\ge 1,\label{eq:5.4}
\end{equation}
{{let its weight be}} $\omega_z(\pi)=(1+z)^t$, {{where t is the number of strict inequalities in \eqn{5.4}. Then for each n and k we have}}
$$
\sum_sb(n,k,s)z^s=\sum_{\pi\in R, \sigma(\pi)=n, \nu(\pi)=k}\omega_z(\pi).
$$
\item[(ii)]
{{Let}} {\bf{P}} {{denote the set of primary partitions. For 
each partition}} $\pi\in${\bf{P}}, 
\begin{equation}
\pi: b_1+b_2+...+b_k, \quad b_i-b_{i+1}\ge 0, \quad \text{for} \quad 
i=1,2,..., k-1, \quad \text{and} \quad b_k\ge k, \label{eq:5.5}
\end{equation}
{{let its weight be}} $\omega_z(\pi)=(1+z)^t$, {{where t is the 
number of strict inequalities in \eqn{5.5}. Then for each n and k we have}}
$$
\sum_sb(n,k,s)z^s=\sum_{\pi\in P, \sigma(\pi)=n, \nu(\pi)=k}\omega_z(\pi).
$$
\end{enumerate}
\end{thm}


\bigskip

\section{A Franklin type involution for Theorem 3}

\medskip
For Euler's celebrated Pentagonal Number Theorem, Fabian Franklin provided 
an elegant combinatorial proof (see \cite{An-book}). Franklin's proof utilized an 
involution on Ferrers graphs of partitions of an integer into distinct 
parts by considering the {\it{base}} and {\it{slope}} on these graphs. 
We will now produce an involution on the Ferrers graphs of basis partitions, 
in which the base and slope are replaced by {\it{bottom block}} and 
{\it{right block}}. 

Given a Ferrers graph of a basis partition, we define a {\it{block}} to be 
the full collection of columns of a given length to the right of 
the Durfee square, or the  full collection of rows of a given length 
below the Durfee square. Thus the part to the right of the Durfee 
square is a collection of blocks, as is the portion below the Durfee 
square. The number of blocks below the Durfee square is equal to the 
number of different parts below the Durfee square, and therefore equal 
to the signature of the basis partition. 

Next we define the length of a block to be its column length if it is 
to the right of the Durfee square, and its row length if it is below 
the Durfee square. Since a basis partition $\pi$ is characterized by the 
property that no column of $\pi_r$ equals a row of $\pi_b$, this is 
equivalent to saying that the blocks in the graph of a basis partition 
all have different lengths. 

Finally define the {\it{right block}} to be the block on the extreme 
right in $\pi_r$, and the {\it{bottom block}} to be the block at the 
very bottom of $\pi_b$. Note that the lengths of the right block and 
the bottom block (if they exist) are different. This leads us to

{\bf{Definition B:}} Given the Ferrers graph of a basis partition, we call 
it a Type B partition, if the bottom block exists, and if its length is 
less than the length of the right block (if the right block exists). 
If the right block does not exist, but the bottom block does, it is 
still a Type B partition. 

{\bf{Definition R:}} Given the Ferrers graph of a basis partition, we call 
it a Type R partition, if the right block exists, and if its length is 
less than the length of the bottom block (if the bottom  block exists). 
If the bottom block does not exist, but the right  block does, it is 
still a Type R partition. 

Thus every Ferrers graph of a basis partition, unless it is just the 
Durfee square, is either a Type B partition or a Type R partition, but 
not both. 

{\bf{Involution:}} Now we define a mapping as follows. If we have a 
Type B partition, 
we move the bottom block, convert its rows into columns, and place it 
on the extreme right thereby making it the (new) right block. Thus we 
create a Type R partition. Conversely, if we have a Type R partition, 
we move the right block, convert its columns into rows, and place it 
at the very bottom of the graph to create a Type B partition. Notice 
that under this mapping, the {\it{parity}} of the number of parts below 
the Durfee square, namely the {\it{parity of the signature}}, changes. The 
only graph on which this mapping does not apply, is the graph consisting 
of just the Durfee square and so this provides a combinatorial proof of 
Theorem 3. 

{\bf{Conjugation as an involution:}} It is worth noting that if $\pi$ is 
a basis partition, then its (Ferrers) conjugate $\pi^*$is also a basis 
partition. This is because basis partitions $\pi$ are characterized by 
the property that no row of $\pi_b$ equals a column of $\pi_r$, and this 
property is preserved under conjugation. Also, this characterization 
implies that $\pi$ and $\pi^*$ are different basis partitions unless 
the graph of $\pi$ is simply a Durfee square. Thus by matching every basis 
basis partition $\pi$ with its conjugate $\pi^*$, we see the number 
$b(n)$ of basis partitions of $n$ is even unless $n$ is a perfect square. 
This parity property of $b(n)$ follows from Theorem 3 but conjugation 
as an involution does not prove Theorem 3.      
\bigskip  
  
\section{Basis partitions with a prescribed signature}

\medskip    
In this section we study the generating function of basis partitions 
with a prescribed signature, and discuss how they are connected to 
the Rogers-Ramanujan partitions. We will use the standard notation
$$
(a;q)_n=(a)_n=\prod^{n-1}_{j=0}(1-aq^j),
$$
for any complex number $a$ and a non-negative integer $n$, and when $|q|<1$ 
$$
(a)_{\infty}=lim_{n\to\infty}(a)_n=\prod^{\infty}_{j=0}(1-aq^j).
$$

The construction in Section 3 showed that we get basis partitions 
from primary partitions (equivalent to Rogers-Ramanujan partitions) 
by sliding blocks from the right of the Durfee square and placing them 
below the Durfee square after converting the columns of these blocks 
into rows. If there are $t$ blocks in the portion $\pi_r$, and $j$ of 
these have to be moved, then the number of choices is 
$$
\bin{t}{j}.
$$   
Thus from this observation, and by collecting the coefficient of $z^j$ 
in Theorem 4, we get the following refinement of Hirschhorn's theorem:

\medskip
\begin{thm}
\label{thm:5}
{{Let}} {\bf{R}} {{denote the set of Rogers-Ramanujan partitions. 
Suppose a partition}} $\pi\in${\bf{R}}, 
\begin{equation}
\pi: b_1+b_2+...+b_k, \quad b_i-b_{i+1}\ge 2, \quad \text{for} \quad 
i=1,2,..., k-1, \quad \text{and} \quad b_k\ge 1,\label{eq:7.1 }
\end{equation}
{{has t strict inequalities in \eqn{5.4}. For a given j, define the weight}}
$$
\omega_j(\pi)=\bin{t}{j}.
$$ 
{{Then for each n and k we have}}
$$
b(n,k,j)=\sum_{\pi\in R, \sigma(\pi)=n, \nu(\pi)=k}\omega_j(\pi).
$$
\end{thm}

Theorem 5 is what we get combinatorially by looking at the sliding 
of blocks. But we get a {\it{different}} connection between 
basis partitions of a prescribed signature and Rogers-Ramanujan partitions 
by evaluating the coefficient of $z^j$ in \eqn{5.2}. 

Denote by $B(n;j)$ the number of basis partitions of $n$ with signature 
$j$ (note the difference with $b(n;k)$). Thus 
\begin{equation}
\sum^{\infty}_{k=0}\frac{q^{k^2}(-zq)_k}{(q)_k}=
\sum_j z^j\sum_nB(n;j)q^n\label{eq:7.1}
\end{equation}
But then    
$$
\sum^{\infty}_{k=0}\frac{q^{k^2}(-zq)_k}{(q)_k}=
\sum^{\infty}_{k=0}\frac{q^{k^2}}{(q)_k}\sum^k_{j=0}z^jq^{j(j+1)/2}
\qbin{k}{j}
$$
\begin{equation}
=\sum^{\infty}_{j=0}\frac{z^jq^{(j(j+1)/2}}{(q)_j}
\sum^{\infty}_{i=0}\frac{q^{(i+j)^2}}{(q)_i}=
\sum^{\infty}_{j=0}\frac{z^jq^{(3j^2+j)/2}}{(q)_j}
\sum^{\infty}_{i=0}\frac{q^{i^2+2ij}}{(q)_i}.\label{eq:7.2}
\end{equation}
By comparing coefficients of $z^j$ in \eqn{7.1} and \eqn{7.2}, we get the 
generating function of basis partitions with a prescribed signature as
\begin{equation}
\sum^{\infty}_{n=0}B(n;j)q^n=
\frac{q^{(3j^2+j)/2}}{(q)_j}
\sum^{\infty}_{i=0}\frac{q^{i^2+2ij}}{(q)_i}.\label{eq:7.3}
\end{equation}
We will now interpret this combinatorially. 

The term 
$$
\frac{q^{(3j^2+j)/2}}{(q)_j}
$$
is the generating function of partitions into $j$ parts with difference 
$\ge 3$ between the parts and smallest part $\ge 2$. The term 
$$
\frac{q^{i^2+2ij}}{(q)_i}
$$
is the generating function for partitions into $i$ parts with difference 
$\ge 2$ between the parts (namely, the Rogers-Ramanujan partitions), but 
with the condition that the smallest part is $\ge 2j+1$. Thus we have 
the following result:

\begin{thm}
\label{thm:6}
{{The number of basis partitions}} $B(n;j)$ {{of 
an integer n with signature j, is equal to the number of vector partitions}} 
$(\pi_3,\pi_2)$ {{of n, where}} $\pi_3$ {{is a partition into j parts 
differing by}} $\ge 3$ {{and least part}} $\ge 2$, {{and}} $\pi_2$ 
{{is a partition into parts differing by}} $\ge 2$ { {and least part}}
$\ge 2j+1$. 
\end{thm}

It is possible to convert these vector partitions into ordinary partitions, 
but then we need to count the resulting Rogers-Ramanujan partitions with 
weights. This is complicated and so we do not pursue it here. 

\bigskip

\section{Complete basis partitions}

\medskip
One of the fundamental properties of partitions into distinct parts is that 
if they are represented as Ferrers graphs, then the conjugate partitions have 
the property that every integer from 1 up to largest part occurs as a part.  
Motivated by this characterization of partitions into distinct parts, we now 
define a {\it{complete basis partition}} $\pi$ to be a basis partition 
with the property that if $|D(\pi)|=k$, then every integer $j$ from 1 to $k-1$ 
occurs either as a row length below $D(\pi)$ or as a column length to the 
right of $D(\pi)$, but not both because $\pi$ is a basis partition ($D(\pi)$ 
= Durfee square of $\pi$). The reason 
for requiring only row (column) lengths up to $k-1$ will become clear below. 
From the 
construction of basis partitions in Section 4, it follows that the an integer 
$j$ between 1 and $|D(\pi)|$ is missed as either a row length below $D(\pi)$ 
or a column length to the right of $D(\pi)$ if and only if in the successive 
rank vector {\bf{r}}=$(r_1,r_2,...,r_k)$, we have $r_j=r_{j+1}$, where we 
formally set $r_{k+1}=0$. Thus we have the following characterization of 
complete basis partitions:

\medskip
\noindent
{\bf{Characterization:}} {\it{A basis partition with a}} $k\times k$ 
{\it{Durfee square is complete if and only if its successive rank vector}}
{\bf{r}}=$(r_1,r_2,...,r_k)$ {\it{has the property that}}
\begin{equation}
r_j\ne r_{j+1}, \quad \text{for} \quad j=1,2,...,k-1.\label{eq:8.1}
\end{equation}

We note that a row (column) of length $k$ is missing precisely when $r_k=0$.

We will now compute the generating function of complete basis partitions by 
connecting them with partitions into distinct parts. For this we consider 
partitions into distinct parts which are primary, that is in their Ferrers 
graph representation, there is nothing below the Durfee square. So let $\pi$ 
be a primary partition into distinct parts with $|D(\pi)|=k$. This Ferrers 
graph has two properties:

(i) To the right of the the Durfee square there is a right-angled 
isosceles triangle of nodes with the two equal sides of the triangle being 
of length $\ge k-1$. This guarantees that there are at least $(3k^2-k)/2$ nodes 
in the Ferrers graph. 

(ii) Every integer between 1 and $k-1$ occurs as a column length to the right 
of $D(\pi)$. 

In \cite{Al0097} we discussed the construction of all partitions into distinct parts 
from primary partitions into distinct parts by means of the sliding operation. 
In doing so we noted in \cite{Al0097} that at most one column of a given length could be 
slid down because the resulting partition has to have distinct parts.  
Here we shall focus of generating all complete basis partitions from primary 
partitions. What this means is that once we slide a column of a certain length 
$j$ to the right of the Durfee square, we have to slide all columns of length 
$j$. Thus for a column of length $j$ on the right, there are precisely two 
choices - to slide the entire set of such columns or not to slide at all. Since 
these choices are independent, and every integer $j$ from 1 to $k-1$ occurs as 
a column length, each primary partition into distinct parts generates at least 
$2^{k-1}$ complete basis partitions. If the graph of $\pi$ has a column of 
length $k$, this will provide two more independent choices. Thus if $b_c(n)$ 
denotes the number of complete basis partitions of $n$, its generating 
function is
\begin{equation}
\sum^{\infty}_{n=0}b_c(n)q^n=
1+\sum^{\infty}_{k=1}\frac{2^{k-1}q^{(3k^2-k)/2}(1+q^k)}{(1-q^k)}.\label{eq:8.2}
\end{equation}

The sliding operation that yielded the complete basis partitions from primary 
basis partitions shows that we could count complete basis basis partitions 
by keeping track of their signature. More precisely, if $b_c(n;j)$ denotes 
the number of complete basis partitions with signature $j$, then 
\begin{equation}
\sum^{\infty}_{n=0}\sum^{\infty}_{j=0}b_c(n;j)z^jq^n=
1+\sum^{\infty}_{k=1}\frac{(1+z)^{k-1}q^{(3k^2-k)/2}(1+zq^k)}{(1-q^k)}.\label{eq:8.3}
\end{equation}
Notice that when $z=-1$, the series in \eqn{8.3} reduces to $1+q$. Thus we have, 

\begin{thm}
\label{thm:7}
{{For each integer}} $n\ge 2$, {{the number of 
complete basis partitions of even signature is equal to the number of 
complete basis partitions of odd signature. In particular, for}} $n\ge 2, 
b_c(n)$ {{is always even}}.
\end{thm}

Theorem 7 can also be proved using the Franklin type involution discussed in 
Section 6. 

The hook lengths of primary partitions into distinct parts yield partitions 
into parts that differ by at least 3, and conversely every partition 
$\tilde\pi$ into parts that differ by at least 3 can be realized in terms of 
hooks of a primary partition $\pi$ into distinct parts. If $|D(\pi)|=k$, 
then $\tilde\pi$ has exactly $k$ parts. With this identification, we have 

\begin{thm}
\label{thm:8}
{{Let}} $D_3$ {{denote the set of partitions whose 
parts differ by at least 3. For a partition}} $\tilde\pi\in D_3$, {{let 
its weight be}} $w(\tilde\pi)=2^{\nu(\tilde\pi) -1}$
{{if}} $\tilde\pi$ {{has 1 as a part, and}} 
$w(\tilde\pi)=2^{\nu(\tilde\pi)}$ {{otherwise, where}} $\nu(\tilde\pi)$ 
{{is the number of parts of}} $\tilde\pi$.  {{Then}}
$$   
b_c(n)=\sum_{\tilde\pi\in D_3, \sigma(\tilde\pi)=n}w(\tilde\pi).
$$
\end{thm}

Theorem 8 is the combinatorial interpretation of \eqn{8.2}. 
  
{\bf{Remark:}} Using \eqn{8.2} or Theorem 8, one can show that 
\begin{equation}
b_c(n)\equiv 0(mod\, 4), \quad \text{for} \quad n\ge 5.\label{eq:8.4}
\end{equation}
More generally, using \eqn{8.2} or Theorem 8 one can study the values of $n$ for 
which $b_c(n)$ is a multiple of $2^k$.
   
\bigskip

\section{Basis partitions among partitions with non-repeating odd 
parts}

\medskip
We conclude the paper with a study of basis partitions among the set $P_{o,d}$ 
of partitions with non-repeating odd parts. There are similarities with the 
theory of basis partitions for unrestricted partitions, but there are very 
interesting differences, and the results quite elegant. That is what motivated 
this study. 

Partitions with non-repeating even parts have been studied extensively owing 
to the famous Lebesgue identity. Comparatively less is known about partitions 
with non-repeating odd parts. Our interest in partitions with non-repeating 
odd parts is due to the fact that if they are represented in terms of 2-modular 
Ferrers graphs (that is Ferrers graphs in which there is a 2 at every node 
except possibly at the last node on the right which would be 1 if the part 
is odd), then the ones can occur only in the corners. Also the conjugate of 
such a graph will also be a 2-modular graph of a partition with non-repeating 
odd parts. By exploiting these properties, we established in \cite{Al2010} a Lebesgue 
type expansion for the two parameter generating function 
$$
\prod^{\infty}_{m=1}\frac{(1+bzq^{2m-1})}{(1-zq^{2m})}
$$
of partitions with non-repeating odd parts and studied connections betwen this 
Lebesgue type identity and various fundamental $q$-hypergeometric identities. 
More recently in \cite{Al2016} we studied the combinatorial properties of partitions in 
$P_{o,d}$ by  means of this 2-modular representation. The results in \cite{Al2010} and \cite{Al2016} 
constitute the first systematic study of partitions in $P_{o,d}$ by means of 
their 2-modular graphs. In this section we will study {\it{minimal basis partitions}}
and {\it{basis partitions}} among
the members of $P_{o,d}$. Results on such partitions in $P_{o,d}$ were reported
in \cite{Al2016}, but here we discuss these two types of partitions in $P_{o,d}$ in greater detail, and also provide a combinatorial method to construct them. Indeed, in \cite{Al2016}, reference is made to this paper for details and the combinatorial construction given below. 

{\it{From now on by a partition we mean one in}} $P_{o,d}$ {\it{and by a Ferrers graph we mean a 2-modular Ferrers graph. We shall from now on refer to the Ferrers graph of unrestricted partitions as {\underbar{ordinary}} Ferrers graphs.}}

Given a 2-modular Ferrers graph of a partition $\pi\in P_{o,d}$, we define 
the {\it{length}} of a row (column) to be the number of nodes in it, and the 
{\it{size}} of a row (column) to be the sum of the nodes in it. Thus the 
row 2,2,1 will have length 3 and size 5. We do not need this distinction 
for ordinary Ferrers graphs. 

Next if $\pi\in P_{o,d}$ has a $k\times k$ Durfee square $D(\pi)$ in its 
2-modular graph, we define its successive rank vector to be 
{\bf{r}}=$(r_1,r_2,...,r_k)$ where 
\begin{equation}
r_j=\,\text{size of the j-th row} \quad - \quad \text{size of the j-th column}.
\label{eq:9.1}
\end{equation}
By a {\it{minimal basis partition}} we mean a partition $\mu\in P_{o,d}$ for 
which the sum of its parts $\sigma(\mu)$ is minimal with respect to a given 
successive rank vector. We will now establish 

\begin{thm}
\label{thm:9}
{{Given any vector of integers}} 
{\bf{r}}=$(r_1,r_2,...,r_k)$, {{there exists a unique minimal basis 
partition}} $\mu$ {{for which}} {\bf{r}} {{is the successive 
rank vector}}. 
\end{thm}

We will establish Theorem 9 by actually constructing $\mu$ from a given 
{\bf{r}}.

Before proceeding with the construction, we note that the Durfee square 
$D(\pi)$ of any partition $\pi\in P_{o,d}$ will have all entries as twos, 
except possibly at the right hand bottom corner where there could be a 1. 
If there is a 1 in the Durfee square, then in the graph there will be no nodes 
to its right in the same row and no nodes below it in the same column.  

We will illustrate the construction of $\mu$ in a sequence of steps with the 
successive rank vector $\mu=(1,-4,1,2,-5)$, but the construction applies in 
general as in Section 4. The contruction, although similar to that in 
Section 4 has an important differences involving parity which we will point 
out below. That is why we give this construction even though there is overlap 
with Section 4. 

\medskip

{\bf{\underbar{Construction of the minimal basis partition:}}}

\medskip

Let the given successive rank vector be {\bf{r}}=$(r_1,r_2,...,r_k)$. 

\noindent
{\bf{Step 1:}} Form a $k\times k$ Durfee square with a 1 at the bottom right 
node and twos everywhere else. This is the minimal basis partition 
corresponding to the vector 
\medskip

\begin{tabbing}
XXXXXXXXXXXXXXXXXXXXXXXX \= YYYYYYYYYYYYYYYYYY\kill
\vbox to 1.1truein{\hbox{${\bold r}=(0,0,0,0,0)$}
                    \vfill} 
\>
\hbox to 1.1truein{\epsfxsize=1.1truein\epsfbox{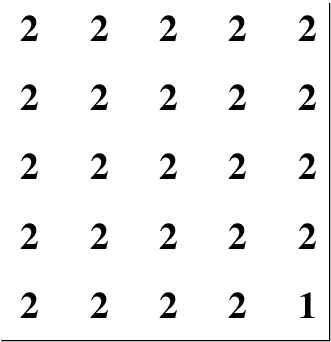}}\\
\end{tabbing}

\noindent
{\bf{Step 2:}} If $r_k=0$, retain the 1 in the Durfee square and move to 
$r_{k-1}$. If $r_k>0$ (resp. $<0$), change the 1 to a 2 in the Durfee square, 
and represent $r_k$ (resp. $|r_k|$)  in 2-modular form as a row (resp. column) 
to the right of (resp. below) the k-th row of the Durfee square. Fill every 
thing above (resp. to the left of) the newly formed row (resp. column) 
with twos. The resulting partition will have the $k$-th successive rank 
as $r_k$, and all other successive ranks will be equal to $r_k$ if $r_k$ is 
even; if $r_k$ is odd, then since the end of the $k$-th row (resp. column) 
will have a 1 with twos elsewhere, all other successive ranks will be $r_k+1$ 
if $r_k>0$ or $r_k-1$ if $r_k<0$. The important thing is that we now have $r_k$ 
as the $k$-th successive rank. We will represent this new vector as 
$(r^*_{k}, r^*_k,...,r^*_k, r_k)$. This is an important difference between 
the situation here and what we had in Section 4. In our example we have 
constructed the minimal basis partition corresponding to 
\medskip

\begin{tabbing}
XXXXXXXXXXXXXXXXXXXXXXXX \= YYYYYYYYYYYYYYYYYY\kill
\vbox to 1.8truein{\hbox{${\bold r}=(-6,-6,-6,-6,-5)$}
                    \vfill} 
\>
\hbox to 1.1truein{\epsfxsize=1.1truein\epsfbox{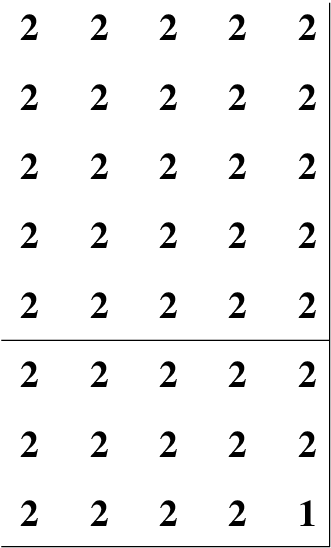}}\\
\end{tabbing}

\noindent
{\bf{Step 3:}} Next consider the difference $r_{k-1}-r^*_k=\delta_{k-1}$. If this 
difference is $>0$, write $\delta_{k-1}$ in 2-modular form to the right of the 
$(k-1)$-st row and fill everything above that with twos. If $\delta_{k-1}<0$, 
write $|\delta_{k-1}|$ in 2-modular form as a column below the 
$(k-1)$-st column 
and fill everything to its left with twos. If $\delta_{k-1}=0$, move on to 
$r_{k-2}$. Thus we have formed the minimal basis partition who last two 
successive ranks are $r_{k-1}$ and $r_k$. The remaining successive ranks 
will all be 
equal, their values being $r_{k-1}$ if it is even, or $r_{k-1}+1$ if $r_{k-1}>0$ 
and odd, or $r_{k-1}-1$ if $r_{k-1}<0$ and odd. We represent these remaining 
successive ranks by $r^*_{k-1}$. Thus we have the minimal basis partition 
corresponding to $(r^*_{k-1},r^*_{k-1},...,r^*_{k-1}, r_{k-1},r_k)$. In our example, 
we have formed the minimal basis partition corresponding to
\medskip

\begin{tabbing}
XXXXXXXXXXXXXXXXXXXXXXXX \= YYYYYYYYYYYYYYYYYY\kill
\vbox to 1.8truein{\hbox{${\bold r}=(2,2,2,2,-5)$}
                    \vfill} 
\>
\hbox to 2.05truein{\epsfxsize=2.05truein\epsfbox{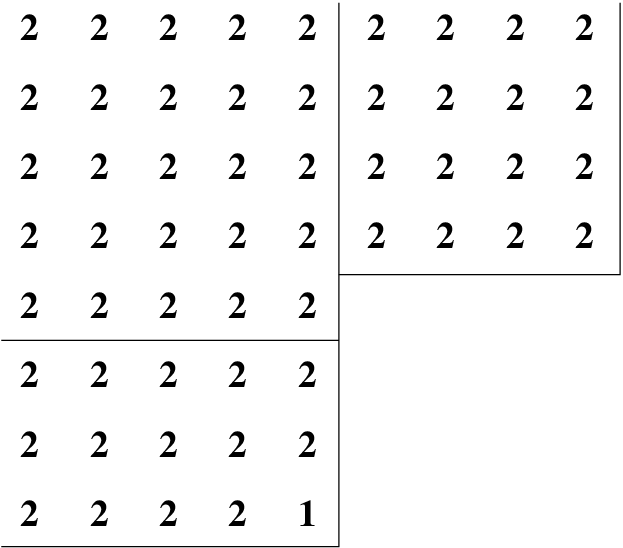}}\\
\end{tabbing}

\noindent
{\bf{Step 4 (General):}}  Consider the difference $r_j-r^*_{j+1}=\delta_i$. 
If this difference is $>0$, represent $\delta_j$ in 2-modular form as a row 
to the right of the $j$-th row, and fill everything above with twos. If the 
difference is $<0$, represent $|\delta_j|$ in 2-modular form as a column 
below the $j$-th column, and fill everything to its left with twos.   
If $\delta_j=0$, proceed to $r_{j-1}$. Note that all successive ranks from 
$r_j$ on are $r_j,r_{j+1},..., r_k$. the remainin g $j-1$ successive ranks are 
all equal and we represent them by $r^*_j$. So we have formed the minimal 
basis partition corresponding to $(r^*_j,r^*_j,...,r^*_j,r_j,r_{j+1},...,r_k)$. 

\medskip
\noindent
{\bf{Step 5:}} Complete the construction by considering the final case $j=1$ 
in Step 4. we then have the minimal basis partition $\mu$ corresponding to the 
given successive rank vector {\bf{r}}. In our example, we have the minimal 
basis partition corresponding to
\medskip

\begin{tabbing}
XXXXXXXXXXXXXXXXXXXXXXXX \= YYYYYYYYYYYYYYYYYY\kill
\vbox to 2.45truein{\hbox{${\bold r}=(1,-4,1,2,-5)$}
                    \vfill} 
\>
\hbox to 2.7truein{\epsfxsize=2.7truein\epsfbox{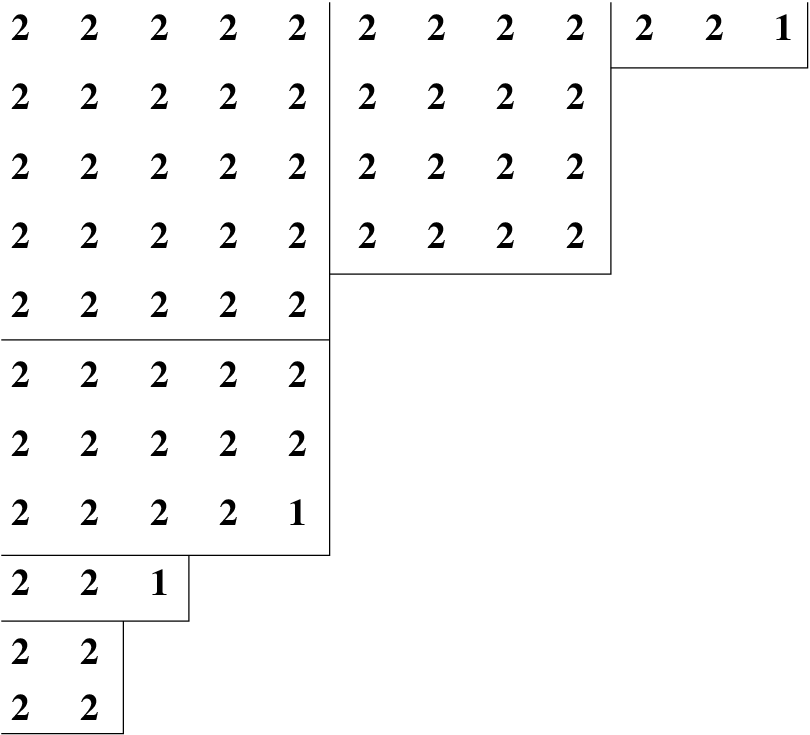}}\\
\end{tabbing}

Our construction has shown that minimal basis partitions are  
characterized by 

\begin{thm}
\label{thm:10}
$\mu$ {is a mimimal basis partition if and only if, 
no row below the Durfee square is equal in {\underbar{length}} 
to a column to the right of the Durfee square, and the right hand 
bottom entry in the Durfee square is 1 if the last successive rank is 0.}
\end{thm}

\medskip

{\bf{Remarks:}}

(1) It is very interesting that even though minimal basis 
partitions were defined in terms of {\it{size}}, the characterization is 
in terms of {\it{length}}. In view of this, we now define a 
{\it{basis partition}} to be a partition in $P_{o.d}$ for which no row 
below the Durfee square is equal \textit{in size} to a column to the right of 
the Durfee square. For example in a basis partition we could have 
2,2,1 as a row and 2,2,2 as a column, but this is not permissible in a minimal 
basis partition. Thus every minimal basis partition is a basis partition, but 
the converse is not true.        

(2) Berkovich and Garvan \cite{Be-Ga02} have studied successive ranks among members of $P_{o,d}$, but they defined the rank as the difference between the lengths of the row and the column, not in terms of the diffrences in sizes as we have considered here. But Berkovich-Garvan did not study basis partitions with their definition.  

\medskip 

{\bf\underbar{Construction of basis partitions from minimal basis 
partitions:}}

\medskip 

The set of all basis partitions can be constructed from the set of all 
minimal basis partitions as we show now. 

Suppose we are given a minimal basis partition $\mu$. If the last successive 
rank is 0, then in the Durfee square of $\mu$ we have a 1. If we replace 
this 1 by a 2, the last rank would still be 0 and so we would get a basis 
partition corresponding to the same successive rank vector. 

Note that in any hook of $\mu$, there can be at most one 1 - either at the 
northeastern end or at the southwestern end but not both. Also, if there is 
a 1 in a hook of $\mu$, then the 2 at the other end cannot be a corner, because 
it is were, then $\mu$ would have a row and column of equal length. So if 
there is a 1 in a hook of $\mu$, we could replace that by a 2, and add a 1 
to the other end of the hook {\underbar{provided}} there is a 2 immediately 
to the left of (or above) the position where are introducing the new 1. We call 
such a hook as an allowable hook. If a minimal basis partition has $\alpha$ 
allowable hooks (including the special case where the Durfee square has a 1 
in it), then with each such hook we have the choice of making such a change 
or not. Thus each minimal basis partition with $\alpha$ allowable hooks, will 
spawn $2^{\alpha}$ basis partitions, including $\mu$, all with the same 
successive rank vector. This procedure will generate all basis 
partitions from minimal basis partitions. 

\medskip

{\bf \underbar{Constructing basis partitions by the sliding operation:}}

\medskip

There is a nice way to get all basis partitions from primary partitions in 
$P_{o,d}$ using the sliding operation. 

In $P_{o,d}$ we define a primary partition to be one such that in its 2-modular 
Ferrers graph there are no nodes below the Durfee square. In \cite{Al2016} the 
construction of all partitions in $P_{o,d}$ from primary partitions is 
discussed. This involves the sliding operation along with the creation of 
secondary partitions generated by primary partitions. For the purpose 
of constructing basis partitions, only primary partitions and the sliding 
operation are needed and we discuss this now. 

Given a primary partition, we consider a column on length $j$ only of twos 
to the right of the Durfee square.  
If there are several such columns of the same length $j$, we either can slide 
all the columns or none, thereby giving us two choices for each collection 
of columns of twos of a given length. If we have a column of length $j$ but 
whose 
size is odd, then such a column can also be slid down giving us once again 
two choices. Through such sliding operations, all basis partitions can be 
generated. Since the weights now will be powers of 2, there is an analogue to 
Hirschhorn's Theorem 1 of section 2, 
but to determine these weights we need to consider seperately columns
having only twos, and columns of twos ending in a one.
But before stating this, we will 
consider the hook operation on primary partitions. 

Given a primary partition in $P_{o,d}$, consider the partition $\tilde\pi$ 
generated by summing the entries on the hooks of $\pi$. This will be a 
partition whose parts will differ by $\ge 4$ with strict inequality if a part 
is odd. We call such partitions as {\it{special}}, 
and denote the set of all special partitions by $\Psi$.
Conversely given a special 
partition, we can find a primary partition whose hook lengths will correspond 
to the given special partition. Thus the hook operation 
\begin{equation}
\pi\mapsto \tilde\pi,\label{eq:9.2}
\end{equation}
is a bijection between primary partitions and special partitions. It is to be 
noted that the existence of a column of length $j$ to the right of the Durfee 
square of a primary partition $\pi$ is equivalent to saying 
that that the difference between the $j$-th part and the $j+1$-st part 
of $\tilde\pi$ is $>4$. Here we adopt the convention that if $\tilde\pi: 
b_1+b_2+...+b_k$, then $b_{k+1}=-2$.  
The construction of basis partitions from primary partitions in the preceding 
paragraph together with \eqn{9.2} yields the following analogue to 
Hirschhorn's theorem:

\begin{thm}
\label{thm:11}
If $\tilde\pi:b_1+b_2+...b_k$ is a special partition, then let
$$
b_i - b_{i+1} \ge  m(b_{i+1}) := 4 + [b_{i+1}]_2,
$$
where $[n]_2=0$ if $n$ is even, and $=1$ if $n$ is odd.
Define the weight $w(\tilde\pi)=2^\ell$, where 
\begin{align*}
\ell &= \{\mbox{number of gaps $ b_i - b_{i+1} >  m(b_{i+1}$)}\} \\
& \qquad + 
\{\mbox{number of gaps $ b_i - b_{i+1} >  m(b_{i+1})+2$ when $b_i$ is odd}\},
\end{align*}
with the convention $b_{k+1}=-2$. 
If $b(n)$ is the number of basis partitions of 
n in $P_{o,d}$, then
\begin{equation}
b(n)=\sum_{\tilde\pi\in \Psi,\sigma(\tilde\pi)=n} w(\tilde\pi).\label{eq:9.3}
\end{equation}
More generally, if $b(n,j)$ is the number of basis partitions in 
$P_{o,d}$ {{with signature j, then}}  
\begin{equation}
\sum_jb(n;j)z^j=\sum_{\tilde\pi\in \Psi, \sigma(\tilde\pi)=n}w_z(\tilde\pi),\label{eq:9.4}
\end{equation}
where $w_z(\tilde\pi)=(1+z)^{\ell}$, if $\tilde\pi$ with $\ell$ 
as above. 
\end{thm}

\medskip     

{\bf\underbar{The generating function of basis partitions and minimal basis 
partitions:}}

\medskip

We will begin with the generating function of basis partitions which is 
easier to derive.  

To get the generating function of basis partitions 
we start with a $k\times k$ Durfee square of all twos which we call Case 2, 
and a $k\times k$ Durfee square of twos with a 1 at the bottom right hand 
corner, which we call Case 1. For Case 2, consider a collection of integers 
all equal to $2j$ represented as a row of twos $j$ in number, 
where $1\le j\le k$. We have 
either a choice of placing ALL of them to the right of the Durfee square as 
columns, or ALL of them as rows below the Durfee square. This gives two choices 
for each set of such $2j$. Similarly if we have an integer $2j-1$ occurring 
singly, we could either place represent it as a column to the right of the 
Durfee square or as a row below the Durfee square. This yields the generating 
function
\begin{equation}
q^{2k^2}\frac{(-q^2;q^2)_k}{(q^2;q^2)_k}.(-2q;q^2)_k.\label{eq:9.5}
\end{equation}
For Case 1 with a $k\times k$ Durfee square, the argument is the same as above 
except that we can only have $1\le j\le k-1$, since $j=k$ is not allowed. This 
yields the generating function
\begin{equation}
q^{2k^2-1}\frac{(-q^2;q^2)_{k-1}}{(q^2;q^2)_{k-1}}.(-2q;q^2)_{k-1}.\label{eq:9.6}
\end{equation}
To get the generating function of $b(n)$, the number of basis partitions of 
$n$, we need to sum the expressions in \eqn{9.5} and \eqn{9.6} over all $k$ to get
\begin{equation}
\sum^{\infty}_{n=0}b(n)q^n=
\sum^{\infty}_{k=0}\frac{q^{2k^2}(-q^2;q^2)_k(-2q;q^2)_k}{(q^2;q^2)_k}
+\sum^{\infty}_{k=1}\frac{q^{2k^2-1}(-q^2;q^2)_{k-1}(-2q;q^2)_{k-1}}{(q^2;q^2)_{k-1}}.
\label{eq:9.7}
\end{equation}

The above combinatorial derivation of the generating function actually allows 
us to refine \eqn{9.7} by introducing parameters $z$ and $b$ whose powers keep 
track of the signature and the the number of odd parts of the partition, 
respectively. That is if $b(n;j,i)$ denotes the number of basis partitions of 
$n$ with signature $j$ and the number of odd parts equal to $i$, then \eqn{9.7} 
refines to
$$
\sum^{\infty}_{k=0}\frac{q^{2k^2}(-zq^2;q^2)_k(-b(1+z)q;q^2)_k}{(q^2;q^2)_k}
+\sum^{\infty}_{k=1}\frac{bq^{2k^2-1}(-zq^2;q^2)_{k-1}(-b(1+z)q;q^2)_{k-1}}
{(q^2;q^2)_{k-1}}
$$
\begin{equation}
=\sum_{n,j,i}b(n;j,i)b^iz^jq^n.\label{eq:9.8}
\end{equation}

To determine the generating function of minimal basis partitions, we need 
to observe that in Case 1, when the Durfee square has dimension $k\times k$, if a collection 
of parts all equal to $2j$ are represented as a set of columns to the 
right of the Durfee square or as a set of rows below the Durfee square, 
but not both, then the integer $2j-1$ if it is to be included in the graph, 
HAS TO BE PLACED ALONGSIDE the collection of $2j$. There is no choice as to 
where to place the $2j-1$ if the $2j$ occurs. But then, if $2j$ DOES NOT OCCUR, 
we could have $2j-1$ represented either as a row below the Durfee square or 
as a column to the right of the square, thereby giving two choices. With 
regard to Case 2, the above observations in this paragraph all hold. In 
addition, we need to note that the last entry in the successive rank vector 
has to be non-zero, thereby forcing the graph to have either a row of 
{\it{length}} $k$ below the Durfee square or a column of {\it{length}} $k$ 
to the right of the square, but not both. This row or column of length $k$ 
could represent either $2k-1$ or $2k$. With these observations, we can 
suitably modify \eqn{9.7} and write down the generating function  of $b_m(n)$, 
the number of minimal basis partitions of $n$, namely:      
$$
\sum^{\infty}_{n=0}b_m(n)q^n=1+ \sum^{\infty}_{k=1}q^{2k^2-1}\prod^{k-1}_{j=1}
\left\{1+2q^{2j-1}+\frac{2q^{2j}(1+q^{2j-1})}{(1-q^{2j})}\right\}
$$  
\begin{equation}
+ \sum^{\infty}_{k=1}q^{2k^2}\prod^{k-1}_{j=1}
\left\{1+2q^{2j-1}+\frac{2q^{2j}(1+q^{2j-1})}{(1-q^{2j})}\right\}
\times
\left(2q^{2k-1}+\frac{2q^{2k}(1+q^{2k-1})}{(1-q^{2k})}\right).\label{eq:9.9}
\end{equation}

As in the case of the generating function of basis partitions, we can refine 
\eqn{9.9}, but here it is best to keep track of the number of different lengths 
below the Durfee square, which we call as $\ell$-{\it{signature}}, which, 
for the sake of clarity, we will keep track by a parameter $\zeta$. Thus we 
have the following refinement of \eqn{9.9}: 
$$
\sum_{n,j}b_m(n;j){\zeta}^jq^n=1+ \sum^{\infty}_{k=1}q^{2k^2-1}\prod^{k-1}_{j=1}
\left\{1+(1+\zeta)q^{2j-1}+\frac{(1+\zeta)q^{2j}(1+q^{2j-1})}{(1-q^{2j})}\right\}+
$$  
\begin{equation}
\sum^{\infty}_{k=1}q^{2k^2}\prod^{k-1}_{j=1}
\left\{1+(1+\zeta)q^{2j-1}+\frac{(1+\zeta)q^{2j}(1+q^{2j-1})}{(1-q^{2j})}\right\}
\times
\left[(1+\zeta)q^{2k-1}+\frac{(1+\zeta)q^{2k}(1+q^{2k-1})}{(1-q^{2k})}\right],\label{eq:9.10}
\end{equation}
where $b_m(n;j)$, the number of minimal basis partitions 
of $n$ with $\ell$-signature equal to $j$.

\medskip

{\bf\underbar{Parity results:}}

\medskip

The series in \eqn{9.8} collapses to 
\begin{equation}
\sum^{\infty}_{k=0}q^{2k^2}+b\sum^{\infty}_{k=1}q^{2k^2-1}.\label{eq:9.11}
\end{equation}
when $z=-1$ and so it yields a nice partity result with a parameter $b$. 
To state this let us denote by $\nu_0(\beta)$ and $\psi(\beta)$ the number 
of odd parts of $\beta$ and the signature of $\beta$, respectively. Here 
$\beta$  is a basis partition in $P_{o,d}$. Also let  
\begin{equation}
B_e(n,b)=\sum_{\beta\in P_{o,d}, \sigma(\beta)=n, \psi(\beta)=\, even}b^{\nu_o(\beta)},\label{eq:9.12}
\end{equation}
and 
\begin{equation}
B_o(n,b)=\sum_{\beta\in P_{o,d}, \sigma(\beta)=n, \psi(\beta)=\, odd}b^{\nu_o(\beta)},\label{eq:9.13}
\end{equation}
where the sums are over basis partitions $\beta$. Then the collapse of \eqn{9.7} 
to \eqn{9.11} has the following partition interpretation: 

\begin{thm}
\label{thm:12}
$$
B_e(n,b)-B_o(n,b)= 1 \,\, \text{if} \,\, n=2k^2, \quad b \,\, \text{if} \,\, 
n=2k^2-1, \quad \, \text{and 0 otherwise.}
$$
\end{thm}

Finally we note that \eqn{9.10} collapses to 
\begin{equation}
1+\sum^{\infty}_{k=1}q^{2k^2-1},\label{eq:9.14}
\end{equation}
when $\zeta=-1$. To interpret this collapse, we denote by $\lambda(\mu)$ 
the $\ell$-signature of a minimal basis partition $\mu$, and by 
\begin{equation}
M_e(n)=\sum_{\mu\in P_{o,d}, \sigma(\mu)=n, \lambda(\mu)=\, even}1,\label{eq:9.15}
\end{equation}
and 
\begin{equation}
M_o(n)=\sum_{\mu\in P_{o,d}, \sigma(\mu)=n, \lambda(\mu)=\, odd}1,\label{eq:9.16}
\end{equation}
where the sum is over minimal basis partitions $\mu$. Then the collapse of 
\eqn{9.10} to \eqn{9.14} has the following interpretation: 

\begin{thm}
\label{thm:13}
$$
M_e(n)-M_o(n)= 1 \,\, \text{if} \,\,  
n=2k^2-1, \quad \, \text{and 0 otherwise.}
$$
\end{thm}

\medskip

{\bf Concluding Remark:} We would like to point out that 
Andrews \cite{An2015}, motivated by the generating function \eqn{2.3} 
of basis partitions, has investigated the more general series
\begin{equation}
G(a,z;q):=\sum^{\infty}_{k=0}\frac{a^kq^{k^2}(-zq)_k}{(q)_k},\label{eq:9.17}
\end{equation}
which gives the generating function of Rogers-Ramanujan partitions in 
the special case $z=0$. Using $G(a,z;q)$ he studied certain polynomials 
which he called {\it{Basis Partition Polynomials}}, and employed these 
to obtain a more rapidly convergent generating function for basis partitions 
as well as a new proof of the Rogers-Ramanujan identities. Note that we 
have combinatorially interpreted $G(1,z;q)$ in Section 7 as the generating 
function of basis partitions where the power of $z$ keeps track of the 
signature. It would be worthwhile to study connections between the 
combinatorics of basis partitions established in this paper and the results 
obtained by Andrews.

\bigskip

{\bf Acknowledgements: } This work which was was begun in 2007 was completed 
only in Fall 2013 when I visited The Pennsylvania State University, at 
which time the results of Section 9 were established. I would like to thank 
George Andrews for arranging that visit and for urging me to write this 
manuscript which was long overdue. But this paper remained unpublished 
since 2013. I thank Gaurav Bhatnagar, Frank Garvan, Christian Krattenthaler, 
and Michael Schlosser, for including this paper in the special volume 
of SIGMA in honor of Steve Milne's 75th birthday. 
Finally, I thank the two referees for their careful reading of the manuscript
and for their helpful suggestions.

\end{document}